\documentclass[a4paper,12pt,twoside]{article}
\usepackage{pst-fill,pst-grad}
\usepackage{textcomp}
\usepackage[english]{babel}
\usepackage{graphicx}
\usepackage{amsmath}
\usepackage{float}
\usepackage[matrix,arrow,curve]{xy}
\usepackage{pstricks}
\usepackage{amsmath,amsfonts,verbatim,afterpage,theorem,euscript,mathrsfs,amssymb}
\usepackage{amsfonts}
\usepackage{amssymb}
\usepackage{array}
\usepackage{dsfont}
\usepackage{hyperref}
\newtheorem{Theoreme}{Theorem}

\newtheorem{Lemme}{Lemma}[section]

\newtheorem{Remarque}{Remark}[section]
\newcommand{\mysection}{\setcounter{equation}{0} \section}

\title{\bf Mixed Sobolev-like Inequalities in Lebesgue spaces of variable exponents and in Orlicz spaces} 
\author{Diego Chamorro\footnote{Laboratoire de Math\'ematiques et Mod\'elisation d'Evry (LaMME) - UMR 8071. Universit\'e d'Evry Val d'Essonne, 23 Boulevard de France, 91037 Evry Cedex, France. email: \textit{diego.chamorro@univ-evry.fr}}}

\begin{document}
\maketitle
\begin{scriptsize}
\abstract{In this short article we show a particular version of the Hedberg inequality which can be used to derive, in a very simple manner, functional inequalities involving Sobolev and Besov spaces in the general setting of Lebesgue spaces of variable exponents and in the framework of Orlicz spaces. }\\[3mm]
\textbf{Keywords:} Hedberg inequality; Sobolev spaces; Besov spaces.\\
\textbf{Mathematics Subject Classification} 46E35 ; 26D10 ; 46E30
\end{scriptsize}
\section{Introduction and presentation of the results}
We study in this article simple proofs for a family of Sobolev-like inequalities using as base spaces the Lebesgue spaces of variable exponent and the Orlicz spaces.\\

Let us start recalling that, for a smooth function $f:\mathbb{R}^n\longrightarrow \mathbb{R}$ with $n\geq 1$ and for a positive parameter $s>0$, we can define the action of the fractional power of the Laplace operator $(-\Delta)^{\frac s2}$ over the function $f$ in the Fourier level by the expression 
\begin{equation}\label{Def_FracLaplace}
\widehat{(-\Delta)^{\frac s2}(f)}(\xi):=|\xi|^s\widehat{f}(\xi),
\end{equation}
and this definition can be extended to tempered distributions $f\in \mathcal{S}'(\mathbb{R}^n)$ by standard procedures (see \cite[Chapter 6]{Grafakos}). For $0<s<+\infty$ and $1<p<+\infty$, we can consider the Sobolev homogeneous space $\dot{W}^{s,p}(\mathbb{R}^n)$ as the subset of $\mathcal{S}'(\mathbb{R}^n)$ such that the quantity $\|f\|_{\dot{W}^{s,p}(\mathbb{R}^n)}:=\|(-\Delta)^\frac{s}{2}(f)\|_{L^p(\mathbb{R}^n)}$ is finite (a more precise definition of these spaces is given in Section \ref{Secc_LebesgueOrlicz} below). \\

Now, if $1<p<+\infty$ and $0<s< n/p$  we have the following classical Sobolev inequality
\begin{equation}\label{Sobolev1}
\|f\|_{L^q(\mathbb{R}^n)}\leq C\|f\|_{\dot{W}^{s,p}(\mathbb{R}^n)},
\end{equation}
where the parameter $q$ is linked to the parameters $s, p$ and to the dimension $n$  by the relationship 
\begin{equation}\label{Sobolev2}
q=\frac{np}{n-sp}\qquad\qquad \mbox{(which can be rewritten as $\frac{1}{q}=\frac{1}{p}-\frac{s}{n}$).}\hspace{-2cm}
\end{equation}
\vspace{-5mm}
\begin{Remarque}
The condition (\ref{Sobolev2}) above can be easily deduced from inequality (\ref{Sobolev1}) by homogeneity with respect to dilations: indeed consider the function $f_\lambda(x)=f(\lambda x)$ with $\lambda>0$, by a change of variables one obtains $\|f_\lambda\|_{L^q(\mathbb{R}^n)}=\lambda^{-\frac nq}\|f\|_{L^q(\mathbb{R}^n)}$ and $\|f_\lambda\|_{\dot{W}^{s,p}(\mathbb{R}^n)}=\lambda^{s-\frac np}\|f\|_{\dot{W}^{s,p}(\mathbb{R}^n)}$. Thus, in order to obtain the estimate (\ref{Sobolev1}) for all $\lambda>0$, we must have the restriction $-\frac nq=s-\frac np$, which is exactly (\ref{Sobolev2}).
\end{Remarque}
Many proofs of this inequality (\ref{Sobolev1}) are available in the litterature and this type of inequalities has been studied in many different settings and admits several generalizations (see \emph{e.g.} \cite{Brezis}, \cite{GMO}).\\

It is classical to link the previous estimate to the Hardy-Littlewood-Sobolev inequality: for $1<p<+\infty$, $0<s< n/p$ and if the parameter $q$ is given by the relationship (\ref{Sobolev2}) above, then we have
\begin{equation}\label{Sobolev3}
\|I_s(f)\|_{L^q(\mathbb{R}^n)}\leq C\|f\|_{L^p(\mathbb{R}^n)},
\end{equation}
where the operator $I_s$ is the Riesz potential defined in the Fourier level by the formula 
\begin{equation}\label{Def_Riesz}
\widehat{I_s(f)}(\xi):=|\xi|^{-s}\widehat{f}(\xi), 
\end{equation}
which is valid for all $s>0$ such that $0<s<n$. See \cite[Chapter 6]{Grafakos} for more details.
\begin{Remarque}
Let us note that the Fourier representation of the fractional powers of the Laplace operator $(-\Delta)^{\frac s2}$ given in (\ref{Def_FracLaplace}) and of the Riesz potential $I_s$ given in (\ref{Def_Riesz}) provides the following semi-group property for $0<s_0, s_1$ and for $0<s_2, s_3<n$ such that $s_2+s_3<n$:
\begin{equation}\label{SemiGroup}
(-\Delta)^{\frac{s_0}{2}}[(-\Delta)^{\frac{s_1}{2}}(f)]=(-\Delta)^{\frac{s_0+s_1}{2}}(f)\quad \mbox{and}\quad I_{s_2}[I_{s_3}(f)]=I_{s_2+s_3}(f).
\end{equation}
Moreover we have the identities
\begin{equation}\label{FracLaplaceRiesz}
I_{s_0}[(-\Delta)^{\frac{s_0}{2}}(f)]=(-\Delta)^{\frac{s_0}{2}}[I_{s_0}(f)]=f,\qquad \mbox{for } 0<s_0<n,
\end{equation}
and
\begin{equation}\label{FracLaplaceRiesz1}
(-\Delta)^{\frac{s_0}{2}}[I_{s_1}(f)]=(-\Delta)^{\frac{s_0-s_1}{2}}(f),\qquad \mbox{for } 0<s_1<s_0<n.
\end{equation}
All these formulas must be taken in the sense of tempered distributions. See also \cite[Exercice 6.1.1]{Grafakos}.
\end{Remarque}

With this remark at hand, the relationship between the Hardy-Littlewood-Sobolev inequality and the classical Sobolev inequality is straightforward: it is enough to consider the function $g=(-\Delta)^{\frac s2}(f)$ in the previous estimate (\ref{Sobolev3}) and to use the definition of the operators $(-\Delta)^{\frac s2}$ and $I_s$ given in 
(\ref{Def_FracLaplace}) and (\ref{Def_Riesz}), as well as the identities (\ref{SemiGroup})-(\ref{FracLaplaceRiesz1}) in order to obtain the inequality (\ref{Sobolev1}).\\

We recall now that the Riesz potentials $I_s$ with $0<s<n$ can also be defined in the real variable by the following expression
\begin{equation*}
I_s(f)(x)=C(n,s)\int_{\mathbb{R}^n}\frac{f(y)}{|x-y|^{n-s}}dy=K_s\ast f(x),
\end{equation*}
where the convolution kernel $K_s$ is given by the locally integrable function $K_s(x):=\frac{1}{|x|^{n-s}}$ (see the book \cite{Grafakos} for more details on the Riesz potential). It is worth noting there that, once we have this characterization of the Riesz potential $I_s$, we can display two easy and straightforward proofs for the Hardy-Littlewood-Sobolev inequality (\ref{Sobolev3}), indeed, assume that $1<p<+\infty$ and $0<s< n/p$: 
\begin{itemize}
\item[$\bullet$] the first proof relies in the fact that the locally integrable function $K_s(x)=\frac{1}{|x|^{n-s}}$ belongs to the Lorentz space $L^{r,\infty}(\mathbb{R}^n)$ with 
\begin{equation}\label{Lorentz1}
r=\frac{n}{n-s},
\end{equation}
(see Section 1.1.1. of the book \cite{Grafakos}) and thus by the Young-O'Neil convolution inequalities (see Theorem 1.4.24 of \cite{Grafakos})  if we define the parameter $q$ by the relationship 
\begin{equation}\label{Lorentz2}
1+\frac{1}{q}=\frac{1}{r}+\frac{1}{p},
\end{equation} 
we can write
\begin{eqnarray}
\|I_s(f)\|_{L^q(\mathbb{R}^n)}&=&\|K_s\ast f\|_{L^q(\mathbb{R}^n)}\leq \|K_s\|_{L^{r,\infty}(\mathbb{R}^n)}\|f\|_{L^p(\mathbb{R}^n)}\label{ONeilYoung}\\
&\leq &C\|f\|_{L^p(\mathbb{R}^n)},\notag
\end{eqnarray}
which is the Hardy-Littlewood-Sobolev inequality (\ref{Sobolev3}).
Note that with the definition of the parameter $r$ given in (\ref{Lorentz1}) and with the relationship (\ref{Lorentz2}) above, we readily obtain the condition (\ref{Sobolev2}) between the indexes $p,q,s$ and the dimension $n$.
\item[$\bullet$] the second proof uses two ingredients: first, the boundedness of the Hardy-Littlewood maximal function $\mathcal{M}$ on Lebesgue spaces. Indeed, we recall here that for a locally integrable function $f:\mathbb{R}^n\longrightarrow \mathbb{R}$, the Hardy-Littlewood maximal function of $f$ is given by 
\begin{equation}\label{MaximalFunction}
\mathcal{M}(f)(x)=\underset{B \ni x}{\sup } \;\frac{1}{|B|}\int_{B }|f(y)|dy,\quad \mbox{where } B \mbox{ is an open ball of } \mathbb{R}^n,
\end{equation}
and then, for $1<p\leq+\infty$, we have the following boundedness property
\begin{equation}\label{MaximalFunction1}
\|\mathcal{M}(f)\|_{L^p(\mathbb{R}^n)}\leq C\|f\|_{L^p(\mathbb{R}^n)}.
\end{equation}
The second ingredient is given by the Hedberg inequality \cite{Hedberg} which reads as follows (for $0<s< n/p$):
\begin{equation}\label{Hedberg}
|I_{s}(f)(x)|\leq C\mathcal{M}(f)(x)^{1-\frac{s p}{n}}\|f\|_{L^{p}(\mathbb{R}^n)}^{\frac{sp}{n}}.
\end{equation}
Thus, taking the $L^q$-norm on the both sides of this estimate, we obtain
$$\|I_{s}(f)\|_{L^q(\mathbb{R}^n)}\leq C\|\mathcal{M}(f)^{1-\frac{s p}{n}}\|_{L^q(\mathbb{R}^n)}\|f\|_{L^{p}(\mathbb{R}^n)}^{\frac{sp}{n}},$$
since we have
\begin{equation*}
\|\mathcal{M}(f)^{1-\frac{s p}{n}}\|_{L^q(\mathbb{R}^n)}=\|\mathcal{M}(f)\|_{L^{q(1-\frac{s p}{n})}(\mathbb{R}^n)}^{1-\frac{s p}{n}},
\end{equation*}
and noting that by the condition (\ref{Sobolev2}) we have the identity $q(1-\frac{s p}{n})=p$, then we obtain 
$$\|I_{s}(f)\|_{L^q(\mathbb{R}^n)}\leq C\|\mathcal{M}(f)\|_{L^p(\mathbb{R}^n)}^{1-\frac{s p}{n}}\|f\|_{L^{p}(\mathbb{R}^n)}^{\frac{sp}{n}},$$
thus, using the boundedness of the maximal function (\ref{MaximalFunction1}) in the Lebesgue spaces we finally obtain
\begin{eqnarray*}
\|I_{s}(f)\|_{L^q(\mathbb{R}^n)}&\leq&C\|\mathcal{M}(f)\|_{L^p(\mathbb{R}^n)}^{1-\frac{s p}{n}}\|f\|_{L^{p}(\mathbb{R}^n)}^{\frac{sp}{n}}\leq C'\|f\|_{L^p(\mathbb{R}^n)}^{1-\frac{s p}{n}}\|f\|_{L^{p}(\mathbb{R}^n)}^{\frac{sp}{n}}\\
&\leq &C'\|f\|_{L^p(\mathbb{R}^n)},
\end{eqnarray*}
which is the Hardy-Littlewood-Sobolev inequality (\ref{Sobolev3}).\\
\end{itemize}
Each one of these proofs has interesting applications and corollaries, and the aim of this article is to develop a variant of the Hedberg inequality (\ref{Hedberg}) in order to study particular versions of the Sobolev inequalities in two different frameworks: the Lebesgue spaces of variable exponent and the Orlicz spaces.
\subsubsection*{The framework of Lebesgue spaces of variable exponent}
We are interested here to study some functional inequalities in the setting of the \emph{Lebesgue spaces of variable exponents} $L^{p(\cdot)}(\mathbb{R}^n)$ which are defined as follows: first consider a function $p:\mathbb{R}^n\longrightarrow [1,+\infty]$, we will say that $p\in \mathcal{P}(\mathbb{R}^n)$ if $p(\cdot)$ is a measurable function and we define 
\begin{equation*}
p^-=\underset{x\in \mathbb{R}^n}{\mbox{inf ess}} \; \{p(x)\} \qquad \qquad \mbox{and} \qquad \qquad  p^+=\underset{x\in \mathbb{R}^n}{\mbox{sup ess}} \; \{p(x)\}.
\end{equation*}
In order to distinguish between variable and constant exponents, we will always denote exponent functions by $p(\cdot)$, moreover, for the sake of simplicity and to avoid technicalities, we will always assume here that we have
$$1<p^-\leq p^+<+\infty.$$ 
Then, for $f:\mathbb{R}^n\longrightarrow \mathbb{R}$ a measurable function we consider  the modular function $\varrho_{p(\cdot)}$ associated with $p\in \mathcal{P}(\mathbb{R}^n)$ by the expression
\begin{equation*}
\varrho_{p(\cdot)}(f)=\int_{\mathbb{R}^n}|f(x)|^{p(x)}dx.
\end{equation*}
Of course, if the function $p(\cdot)$ is constant we obtain the classical Lebesgue spaces and we can derive from the modular function  $\varrho_p$ a norm defined by 
$$\|f\|_{L^{p}(\mathbb{R}^n)}=\left(\int_{\mathbb{R}^n}|f(x)|^{p}dx\right)^{\frac1p}.$$
As we can easily guess from the fact that $p(\cdot)$ is a function, we cannot simply replace in the previous formula the constant exponent $\frac1p$ outside the integral by $\frac{1}{p(\cdot)}$ and in order to overcome this issue we consider the Luxemburg norm given by the expression
\begin{equation}\label{Def_LuxNormLebesgue}
\|f\|_{L^{p(\cdot)}(\mathbb{R}^n)}=\inf\{\lambda > 0: \, \varrho_{p(\cdot)}(f/\lambda)\leq1\},
\end{equation}
and we will define the spaces $L^{p(\cdot)}(\mathbb{R}^n)$ as the set of measurable functions such that the quantity $\|\cdot\|_{L^{p(\cdot)}(\mathbb{R}^n)}$ is finite.\\ 

Although the spaces $L^{p(\cdot)}(\mathbb{R}^n)$ have some nice structural properties (derived from the fact that they are normed spaces, see the books \cite{Cruz1}, \cite{Diening} for more details), some of the usual tools between these spaces are delicate to use and thus, if we want to study Sobolev or Hardy-Littlewood-Sobolev inequalities in this framework, some slightly different ideas must be used.\\ 

For example, the approach based in convolution given in (\ref{ONeilYoung}) seems hard to display as the Young inequalities for convolution are not available in general, indeed, if $p\in \mathcal{P}(\mathbb{R}^n)$, it is known that the estimate
$$\|f\ast g\|_{L^{p(\cdot)}(\mathbb{R}^n)}\leq C\|f\|_{L^{p(\cdot)}(\mathbb{R}^n)}\|g\|_{L^1(\mathbb{R}^n)},$$ 
holds true for all $f\in L^{p(\cdot)}(\mathbb{R}^n)$ and all $g\in L^1(\mathbb{R}^n)$ if and only if $p(\cdot)$ is constant. See Corollary 3.6.4 of \cite{Diening}, see also Theorem 5.19 and Example 5.21 of the book \cite{Cruz1} for a proof of this fact. However, let us mention for the sake of completeness that some weak versions of this convolution inequality are available: see Proposition 5.20 of \cite{Cruz1}, or Theorem 3.6.5 of \cite{Diening}.\\

On the other hand, a direct adaptation of the Hedberg inequality (\ref{Hedberg}) to this framework will inevitably introduce some variable exponents which have to be studied carefully and make the proof much harder. Moreover, some extra assumptions on the function $p(\cdot)$ are needed to obtain the boundedness of the Hardy-Littlewood maximal operator $\mathcal{M}$ (which is the key argument of this approach). Indeed, we will say that a measurable function $p\in \mathcal{P}(\mathbb{R}^n)$ belongs to the class $\mathcal{P}^{log}(\mathbb{R}^n)$ if we have
$$\left|\frac{1}{p(x)}-\frac{1}{p(y)}\right|\leq \frac{C}{\log(e+1/|x-y|)}\qquad \mbox{for all } x,y\in \mathbb{R}^n,$$
and if 
$$\left|\frac{1}{p(x)}-\frac{1}{p_\infty}\right| \leq \frac{C}{\log(e+|x|)}\qquad \mbox{for all } x\in \mathbb{R}^n,$$
where $\frac{1}{p_\infty}=\underset{|x|\to +\infty}{\lim}\frac{1}{p(x)}$. See Definition 4.1.1 and Definition 4.1.4 of \cite{Diening}.\\

Thus the condition $p\in \mathcal{P}^{log}(\mathbb{R}^n)$ ensures the fact that the Hardy-Littlewood maximal function is bounded in the Lebesgue spaces of variable exponents:
\begin{equation*}
\|\mathcal{M}(f)\|_{L^{p(\cdot)}(\mathbb{R}^n)}\leq C\|f\|_{L^{p(\cdot)}(\mathbb{R}^n)}, \qquad p\in  \mathcal{P}^{log}(\mathbb{R}^n).
\end{equation*}
See Theorem 4.3.8 of \cite{Diening} for a proof of this estimate.\\

With this boundedness property, it is possible to generalize the Hedberg inequality in the setting of the Lebesgue spaces of variable exponents. This strategy to obtain Hardy-Littlewood-Sobolev inequalities was (to the best of our knowledge) first displayed  in \cite{Capone}: if $0<s<n/p^+$ is a parameter and if $p \in \mathcal{P}^{log}(\mathbb{R}^n)$ is a measurable function such that $1<p^-\leq p^+<+\infty$ then we have the inequality 
\begin{equation}\label{PotentialRieszVariable}
\|I_s(f)\|_{L^{q(\cdot)}(\mathbb{R}^n)}\leq C\|f\|_{L^{p(\cdot)}(\mathbb{R}^n)},
\end{equation}
where the function $q(\cdot)$ is given pointwise by the relationship
\begin{equation}\label{PotentialRieszVariable1}
\frac{1}{q(x)}=\frac{1}{p(x)}-\frac{s}{n},
\end{equation}
and we can see that this previous condition is the ``variable exponent'' version of the relationship (\ref{Sobolev2}) stated before. Of course, from this estimate and using standard arguments we can derive the usual Sobolev inequality in the setting of Lebesgue spaces of variable exponents:
\begin{equation*}
\|f\|_{L^{q(\cdot)}(\mathbb{R}^n)}\leq C\|f\|_{\dot{W}^{s,p(\cdot)}(\mathbb{R}^n)},
\end{equation*}
where we have $\|f\|_{\dot{W}^{s,p(\cdot)}(\mathbb{R}^n)}:=\|(-\Delta)^{\frac{s}{2}}(f)\|_{L^{p(\cdot)}(\mathbb{R}^n)}$.\\

Let us remark that in order to obtain (\ref{PotentialRieszVariable}), the authors of \cite{Capone} develop a very interesting theory related to \emph{fractional} maximal functions. In particular, they use some local\footnote{\emph{i.e.} in bounded domains.} Hedberg-like inequalities for this type of operators (see the details in Proposition 3.3 of \cite{Capone}) since unbounded versions of these arguments introduce some problems as mentioned in Example 3.4 of the cited article.\\

To circumvene these problems, another proof of the  Hardy-Littlewood-Sobolev inequalities was given in the book \cite{Diening} and it relies in a different version of the Hedberg estimate. Indeed, under the condition (\ref{PotentialRieszVariable1}) above for the functions $p(\cdot)$, $q(\cdot)$ and with $p \in \mathcal{P}^{log}(\mathbb{R}^n)$, the following inequality (see \cite[ Lemma 6.1.8.]{Diening}) is obtained for functions such that $\|f\|_{L^{p(\cdot)}}\leq 1$:
\begin{equation}\label{HedbergLebesgueVariableDiening}
|I_s(f)(x)|^{q(x)}\leq C\mathcal{M}(f)(x)^{p(x)}+h(x),\quad x\in \mathbb{R}^n,
\end{equation}
where $q(\cdot)$ satisfies (\ref{PotentialRieszVariable1}) and $h\in L^1\cap L^\infty(\mathbb{R}^n)$ is of the following form
$$h(x)\simeq \mathcal{M}\left((e+|x|)^{-n}\right)(x)^{p^-}+(e+|x|)^{-m},\qquad \mbox{with } m>n.$$
As we can see, this version (\ref{HedbergLebesgueVariableDiening}) of the Hedberg inequality in the setting of Lebesgue spaces of variable exponents is slightly different from (\ref{Hedberg}) as a function $h$ appears in the right-hand side. However, and despite of this extra term, with the pointwise inequality (\ref{HedbergLebesgueVariableDiening}) the estimate (\ref{PotentialRieszVariable}) can be easily deduced (see \cite[Theorem 6.1.9.]{Diening}).\\

In this article we will present a simpler version of the Hedberg inequality that will lead us to a mixed variable-constant Sobolev inequality.
\subsubsection*{Presentation of our results}
Our main result gives an alternative version of the Hedberg inequality (\ref{Hedberg}) which is, to our belief, more natural than (\ref{HedbergLebesgueVariableDiening}) in view to obtain generalized Sobolev embeddings:
\begin{Theoreme}[Modified Hedberg Inequality]\label{Theo1}
Consider two positive parameters $s, s_1$ such that $0<s_1<s<n$ and consider a smooth measurable function $f:\mathbb{R}^n\longrightarrow \mathbb{R}$ that belongs to the Besov space $\dot{B}^{-\beta, \infty}_{\infty}(\mathbb{R}^n)$ for some $\beta>0$. Then we have the following version of the Hedberg inequality
\begin{equation}\label{HedbergNuevo1}
|(-\Delta)^{\frac{s_1}{2}}(f)(x)|^{\frac{1}{1-\theta}}\leq C\mathcal{M} \left((-\Delta)^{\frac{s}{2}} (f)\right)(x)\|f\|^{\frac{\theta}{1-\theta}}_{\dot{B}^{-\beta, \infty}_{\infty}(\mathbb{R}^n)},
\end{equation}
where the parameter $0<\theta<1$ is given by the relationship
\begin{equation*}
\frac{s-s_1}{\beta+s}=\theta.
\end{equation*}
\end{Theoreme}
Let us note that, contrary to the usual Hedberg inequality (\ref{Hedberg}) or to the version (\ref{HedbergLebesgueVariableDiening}) with variable exponents given above, all the exponents here depend only on the parameters $s, s_1$ and $\beta$ which are \emph{constants}. Note also that the Lebesgue norm $\|\cdot\|_{L^{p}(\mathbb{R}^n)}$ in (\ref{Hedberg}) is now replaced by a Besov norm $\|\cdot\|_{\dot{B}^{-\beta,\infty}_\infty}$ which will allows us to perform simpler computations. See formula (\ref{Def_BesovThermic}) below for a precise definition of this Besov space.\\

This inequality will lead us, via simple arguments, to the following estimate.
\begin{Theoreme}[Sobolev-like inequalities]\label{Theo2}
Let $p\in \mathcal{P}^{\log}(\mathbb{R}^n)$ with $1<p^-\leq p^+<+\infty$ and let $0<s<n/p^+$. Consider $f:\mathbb{R}^n\longrightarrow \mathbb{R}$ a measurable function such that $f\in \dot{W}^{s,p(\cdot)}(\mathbb{R}^n)$ and $f\in \dot{B}^{-\beta,\infty}_{\infty}(\mathbb{R}^n)$ for some $\beta>0$. For any fixed $s_1\geq 0$ such that $0\leq s_1<s$ we define $\theta=\frac{s-s_1}{\beta+s}<1$ and $q(\cdot)=\frac{p(\cdot)}{1-\theta}$. Then we have $f\in \dot{W}^{s_1,q(\cdot)}(\mathbb{R}^n)$ and the following inequality holds
\begin{equation*}
\|f\|_{\dot{W}^{s_1,q(\cdot)}(\mathbb{R}^n)}\leq C \|f\|_{\dot{W}^{s,p(\cdot)}(\mathbb{R}^n)}^{1-\theta}\|f\|_{\dot{B}^{-\beta,\infty}_{\infty}(\mathbb{R}^n)}^\theta.
\end{equation*}
\end{Theoreme}
This type of inequalities are quite useful in the field of PDEs and can be inserted in the global family of Gagliardo-Nirenberg estimates. In the classical setting of constant coefficients, these estimates are more precise and robust than the usual Sobolev inequalities and they are often known as \emph{improved Sobolev inequalities} or \emph{refined Sobolev inequalites}, see \cite{Bahouri}, \cite{Kolyada}, \cite{GMO} and the references there in for some applications of these estimates.\\

In the setting of Lebesgue spaces of variable exponents,  we will give two different proofs of this inequality: the first proof will use the modified Hedgberg inequality (\ref{HedbergNuevo1}) above while the second proof will use the Littlewood-Paley theory and it will be given in the Appendix \ref{AppendixA}.\\

To the best of our knowledge these results are new in the setting of Lebesgue spaces of variable exponents (and in the setting of Orlicz spaces that will be studied in Section \ref{Secc_Orlicz} below) and we will see that we can derive from these estimates some inequalities that may have their own interest.\\

Indeed, if $p\in \mathcal{P}^{log}(\mathbb{R}^n)$ is a variable exponent with $1<p^-\leq p^+<+\infty$ and $1<\mathfrak{p}<+\infty$ is a constant, we consider the intersection space $\mathcal{L}^{p(\cdot)}_\mathfrak{p}(\mathbb{R}^n):=L^{p(\cdot)}(\mathbb{R}^n)\cap L^{\mathfrak{p}}(\mathbb{R}^n)$, that can be normed by quantity
\begin{equation}\label{MixedLebesgue}
\|\cdot\|_{\mathcal{L}^{p(\cdot)}_\mathfrak{p}(\mathbb{R}^n)}=\max\{\|\cdot\|_{L^{p(\cdot)}(\mathbb{R}^n)}, \|\cdot\|_{L^{\mathfrak{p}}(\mathbb{R}^n)}\},
\end{equation}
of course, if $p(\cdot)=\mathfrak{p}$, we have $\mathcal{L}^{p(\cdot)}_\mathfrak{p}(\mathbb{R}^n)=L^\mathfrak{p}(\mathbb{R}^n)$.\\

In the same spirit, with $p(\cdot)$ and $\mathfrak{p}$ as above and for $0<s<+\infty$ we define the mixed-norm Sobolev space $\dot{\mathcal{S}}^{s, p(\cdot)}_\mathfrak{p}(\mathbb{R}^n):=\dot{W}^{s,p(\cdot)}(\mathbb{R}^n)\cap \dot{W}^{s,\mathfrak{p}}(\mathbb{R}^n)$, which can be characterized by the functional
\begin{equation}\label{MixedSobolev}
\|\cdot\|_{\dot{\mathcal{S}}^{s, p(\cdot)}_\mathfrak{p}(\mathbb{R}^n)}=\max\{\|\cdot\|_{\dot{W}^{s,p(\cdot)}(\mathbb{R}^n)}, \|\cdot\|_{\dot{W}^{s,\mathfrak{p}}(\mathbb{R}^n)}\}.
\end{equation}
With this notation, a simple corollary of the previous theorem is the following one:
\begin{Theoreme}[Mixed Sobolev inequalities]\label{Theo3}
Consider a constant exponent $1<\mathfrak{p}<+\infty$ and a variable exponent $p\in \mathcal{P}^{\log}(\mathbb{R}^n)$ such that $1<p^-\leq p^+<+\infty$ and fix a parameter $0<s<\min\{n/p^+, n/\mathfrak{p}\}$. Assume that $f\in \dot{\mathcal{S}}^{s, p(\cdot)}_\mathfrak{p}(\mathbb{R}^n)$, then we have the inequality
\begin{equation*}
\|f\|_{L^{\sigma(\cdot)}(\mathbb{R}^n)}\leq C \|f\|_{\dot{\mathcal{S}}^{s, p(\cdot)}_\mathfrak{p}(\mathbb{R}^n)},
\end{equation*}
where $\theta=\frac{s\mathfrak{p}}{n}$ and where the function $\sigma(\cdot)$ satisfies the following condition
\begin{equation}\label{Inegalite2Condition}
\sigma(\cdot)=\frac{np(\cdot)}{n-s\mathfrak{p}}.
\end{equation}
\end{Theoreme}
Another corollary of Theorem \ref{Theo2} is the following:
\begin{Theoreme}[Mixed Hardy-Littlewood-Sobolev inequalities]\label{Theo4}
Consider a constant exponent $1<\mathfrak{p}<+\infty$ and a variable exponent $p\in \mathcal{P}^{\log}(\mathbb{R}^n)$ such that $1<p^-\leq p^+<+\infty$ and fix a parameter $0<s<\min\{n/p^+, n/\mathfrak{p}\}$. Assume that $f\in \mathcal{L}^{p(\cdot)}_\mathfrak{p}(\mathbb{R}^n)$, then we have the inequality
\begin{equation}\label{Inegalite3}
\|I_s(f)\|_{L^{\sigma(\cdot)}(\mathbb{R}^n)}\leq C \|f\|_{\mathcal{L}^{p(\cdot)}_\mathfrak{p}(\mathbb{R}^n)},
\end{equation}
where $\theta=\frac{s\mathfrak{p}}{n}$ and the function $\sigma(\cdot)$ satisfies the same condition (\ref{Inegalite2Condition}).
\end{Theoreme}
Let us remark that if we set $p(\cdot)=\mathfrak{p}$ we recover the classical Sobolev inequalities (\ref{Sobolev1}) in the framework of usual Lebesgue spaces. We note also that, since no simple description seems to be available in unbounded domains of the intersection of a classical Lebesgue space with a variable exponent one (here $\mathfrak{p}$ needs not to be related to $p^-$ or $p^+$), we believe that these mixed-norm results can be useful in some applications.

\subsubsection*{Outline of the article}
The plan of the article is the following: in Section \ref{Secc_LebesgueOrlicz} we introduce the Lebesgue spaces of variable exponents and their main properties, in Section \ref{Secc_DemoTheo1} we give a proof for Theorem \ref{Theo1}. In Section \ref{Secc_ProofTheo2} we give the proofs of  Theorem \ref{Theo2} and in Section \ref{Secc_ProofTheo3} we give the proofs of  Theorem \ref{Theo3} and \ref{Theo4}. In Section \ref{Secc_Orlicz} we study these results in the framework of Orlicz spaces and finally in Section \ref{Secc_General} we will give some possible extensions of our work.
\section{Functional spaces of variable exponents}\label{Secc_LebesgueOrlicz}
We give in this section the precise definition of all the functional spaces involved in our theorems. 

\begin{enumerate}
\item[$\bullet$] \textbf{Lebesgue spaces of variable exponents}. We have already seen in the introduction that for a measurable function $p\in \mathcal{P}(\mathbb{R}^n)$ we can define the Lebesgue space $L^{p(\cdot)}(\mathbb{R}^n)$ as the set of measurable such that the Luxemburg norm $\|\cdot\|_{L^{p(\cdot)}(\mathbb{R}^n)}$ given in (\ref{Def_LuxNormLebesgue}) is finite. Here are some useful properties of these spaces:
\begin{itemize}
\item[(i)] As the quantity $\|\cdot\|_{L^{p(\cdot)}(\mathbb{R}^n)}$ is a norm, for all constant $\lambda\in \mathbb{R}$ we obviously have the identity
\begin{equation}\label{ConstanteNormeLuxem1}
\|\lambda f\|_{L^{p(\cdot)}(\mathbb{R}^n)}=|\lambda|\;\| f\|_{L^{p(\cdot)}(\mathbb{R}^n)}.
\end{equation}
\item[(ii)] The Luxemburg norm is order preserving: if $f,g\in L^{p(\cdot)}(\mathbb{R}^n)$ are such that $|f|\leq |g|$ a.e., then we have 
\begin{equation}\label{OrdreNormeLuxem1}
\|f\|_{L^{p(\cdot)}(\mathbb{R}^n)}\leq \|g\|_{L^{p(\cdot)}(\mathbb{R}^n)},
\end{equation}
see \cite[Proposition 2.7]{Cruz1}.
\item[(iii)] Another particular feature of the Luxemburg norm for Lebesgue spaces of variable exponent is the following: for all real parameter $\alpha>0$ such that $\frac{1}{p^-}\leq \alpha<+\infty$, we have the identity
\begin{equation}\label{PuissanceNorme}
\||f|^\alpha\|_{L^{p(\cdot)}(\mathbb{R}^n)}=\|f\|^\alpha_{L^{\alpha p(\cdot)}(\mathbb{R}^n)}.
\end{equation}
See \cite[Proposition 2.18]{Cruz1} for a proof of this fact.
\end{itemize}
For a more detailed study of these spaces see the books \cite{Cruz1} and \cite{Diening}.\\

\item[$\bullet$]\textbf{Sobolev spaces}. For a measurable function $p\in \mathcal{P}(\mathbb{R}^n)$ such that $1<p^-\leq p^+<+\infty$ and for a positive index $0<s<n/p^+$, we define the Sobolev spaces $\dot{W}^{s,p(\cdot)}(\mathbb{R}^n)$ as the closure of smooth functions with respect to the functional $\|\cdot\|_ {\dot{W}^{s,p(\cdot)}}$ which is given by:
\begin{equation}\label{Def_Sobolev}
\|f\|_ {\dot{W}^{s,p(\cdot)}(\mathbb{R}^n)} =\|(-\Delta)^{\frac{s}{2}}(f)\|_{L^{p(\cdot)}(\mathbb{R}^n)}.
\end{equation}
Now, if $k$ is an integer and if we have $\frac{n}{p^+}+k<s<\frac{n}{p^+}+k+1$, then we shall define the homogeneous Sobolev space $\dot{W}^{s,p(\cdot)}(\mathbb{R}^n)$ as the set of tempered distributions, modulo polynomials of degree $k$, such that quantity (\ref{Def_Sobolev}) is finite.

This precaution in the definition of the homogeneous Sobolev spaces is necessary because the functional (\ref{Def_Sobolev}) is only a semi-norm. Remark in particular that we always have the condition $0<s<n/p^+$ in our results. Of course this definition remains the same for Sobolev spaces of constant exponents.

\item[$\bullet$] \textbf{Besov spaces}. There are many different (and equivalent) ways to define these spaces and we will use the following \emph{thermic} characterization. Indeed, we will define the Besov spaces of indices $(-\beta, \infty, \infty)$ that appear in the previous inequalities as the set of tempered distributions such that the quantity
\begin{equation}\label{Def_BesovThermic}
\|f\|_{\dot{B}^{-\beta,\infty}_{\infty}(\mathbb{R}^n)}=\underset{t>0}{\sup}\;\;t^{\frac \beta2 } \|h_{t}\ast f\|_ {L^\infty(\mathbb{R}^n)},
\end{equation}
if finite. In the previous formula $h_t$ is the heat (or gaussian) kernel. See \cite{Grafakos},  \cite{Stein1} or \cite{Triebel} for more details and equivalent characterization of Besov spaces. Among many properties of the functional $\|\cdot\|_{\dot{B}^{-\beta,\infty}_{\infty}}$ defined above, we will use the following property which is valid for any $s>0$:
\begin{equation}\label{EquivBesov}
\|(-\Delta)^{\frac s2}(f)\|_{\dot{B}^{-\beta-s,\infty}_{\infty}(\mathbb{R}^n)}\simeq \|f\|_{\dot{B}^{-\beta,\infty}_{\infty}(\mathbb{R}^n)}.
\end{equation}
\end{enumerate}

\section{Proof of Theorem \ref{Theo1}}\label{Secc_DemoTheo1}
Recall that we have $0<s_1<s$. We start using the Riemann-Liouville characterization of the positive powers of the Laplacian:  for any integer $k>s/2>s_1/2>0$ we can write
\begin{equation*}
(-\Delta)^{\frac {s_1}{2}}(f)(x)=\frac{1}{\Gamma(k-s_1/2)}\int_{0}^{+\infty}t^{k-\frac{s_1}{2}-1}(-\Delta)^k (h_t\ast f)(x)dt,
\end{equation*}
where $h_t(x)=\frac{1}{(4\pi t)^{\frac{n}{2}}}e^{-\frac{|x|^2}{4t}}$ with $t>0$ is the usual gaussian kernel. Now, by introducing a cut-off parameter $T$ that will be defined below we have
\begin{eqnarray}
|(-\Delta)^{\frac{s_1}{2} }(f)(x)|&\leq &\frac{1}{\Gamma(k-s_1/2)}\left(\int_{0}^{T}t^{k-\frac{s_1}{2}-1}|(-\Delta)^k(h_t\ast f)(x)|dt\right.\notag\\
&& \qquad \qquad\qquad\left.+ \int_{T}^{+\infty}t^{k-\frac{s_1}{2}-1}|(-\Delta)^k(h_t\ast f)(x)|dt\right).\label{Hedberg_Inequality_1}
\end{eqnarray}
We will study these two integral separately and first we will use the following classical result:
\begin{Lemme}\label{Lemma_Maximal_1}
Let $f\in L^{1}_{loc}(\mathbb{R}^n)$ and $\varphi\in \mathcal{S}(\mathbb{R}^n)$. We denote by $\mathcal{M}_{\varphi}(f)$ the maximal function of $f$ (with respect to $\varphi$) which is given by the expression
\begin{equation*}
\mathcal{M}_{\varphi}(f)(x)=\underset{t>0}{\sup}\{|f\ast\varphi_{t}(x)|\}, \quad \mbox{with } \varphi_{t}(x)=t^{-n/2}\varphi(t^{-1/2}x).
\end{equation*}
If the function $\varphi$ is such that $|\varphi(x)|\leq C(1+|x|)^{-n-\varepsilon}$ for some $\varepsilon>0$, then we have the following pointwise inequality
\begin{equation*}
\mathcal{M}_{\varphi}(f)(x)\leq C \mathcal{M}(f)(x),
\end{equation*}
where $\mathcal{M}(f)$ is the Hardy-Littlewood maximal function defined by (\ref{MaximalFunction}).
\end{Lemme}
For a proof of this lemma see \cite[Theorem 2.1.10]{Grafakos}.\\

With this lemma in mind we will study the terms inside the first integral of (\ref{Hedberg_Inequality_1}). Indeed, since $k>s/2$, we remark that we have the identity 
$$(-\Delta)^k(h_t\ast f)=(-\Delta)^{k-\frac{s}{2}}(h_t) \ast (-\Delta)^{\frac{s}{2}}(f).$$
Now, by homogeneity we obtain the identity 
$$(-\Delta)^{k-\frac{s}{2}}(h_t)(x)=t^{-k+\frac{s}{2}} \big((-\Delta)^{k-\frac{s}{2}}h_t\big)(x),$$
and if we denote the function $\varphi_t$ by 
$$\varphi_t(x)= \big((-\Delta)^{k-\frac{s}{2}}h_t\big)(x),$$
we have that $\varphi_t(x)=t^{-n/2}\varphi(t^{-1/2}x)$, moreover, since the heat kernel $h_t$ is a smooth function, with the previous notation we obtain the estimate $|\varphi(x)|\leq C(1+|x|)^{-n-\varepsilon}$. Then we can write
$$(-\Delta)^k(h_t\ast f)(x)= t^{-k+\frac{s}{2}} \varphi_t \ast (-\Delta)^{\frac{s}{2}}(f)(x),$$
and applying the Lemma \ref{Lemma_Maximal_1} we have the following pointwise inequality 
\begin{eqnarray}
|(-\Delta)^k(h_t\ast f)(x)|&=&t^{-k+\frac{s}{2}}| \varphi_t \ast (-\Delta)^{\frac{s}{2}}(f)(x)|\notag\\
&\leq&t^{-k+\frac{s}{2}}\underset{t>0}{\sup}\{| \varphi_t \ast (-\Delta)^{\frac{s}{2}}(f)(x)|\}=t^{-k+\frac{s}{2}} \mathcal{M}_\varphi(f)(x)\notag\\
&\leq & Ct^{-k+\frac{s}{2}}\mathcal{M}\left( (-\Delta)^{\frac{s}{2}}(f)\right) (x).\label{Estimation1}
\end{eqnarray}
Next, we study the terms inside the second integral of (\ref{Hedberg_Inequality_1}) and we simply use the fact that
\begin{equation*}
|(-\Delta)^k (h_t\ast f)(x)|=|h_t\ast (-\Delta)^k (f)(x)|\leq C t^{\frac{-\beta-2k}{2}}\|(-\Delta)^k (f)\|_{\dot{B}^{-\beta-2k, \infty}_{\infty}(\mathbb{R}^n)},
\end{equation*}
which is a consequence of the definition (\ref{Def_BesovThermic}) of the Besov spaces $\dot{B}^{-\beta-2k, \infty}_{\infty}(\mathbb{R}^n)$. Then we use the equivalence (\ref{EquivBesov}) to obtain
\begin{equation}\label{Estimation2}
|(-\Delta)^k (h_t\ast f)(x)|\leq C t^{\frac{-\beta-2k}{2}}\|f\|_{\dot{B}^{-\beta, \infty}_{\infty}(\mathbb{R}^n)}.
\end{equation}
With these two inequalities (\ref{Estimation1}) and (\ref{Estimation2}) at hand, we apply them in (\ref{Hedberg_Inequality_1}) and one has
\begin{eqnarray*}
|(-\Delta)^{\frac{s_1}{2}}(f)(x)|&\leq &\frac{C}{\Gamma(k-s_1/2)}\left(\int_{0}^{T} t^{k-\frac{s_1}{2}-1} t^{-k+\frac{s}{2}}\mathcal{M}\left( (-\Delta)^{\frac{s}{2}}(f)\right) (x)dt\right.\\
&&\left.+ \int_{T}^{+\infty}t^{k-\frac{s_1}{2}-1} t^{\frac{-\beta-2k}{2}}\|f\|_{\dot{B}^{-\beta, \infty}_{\infty}(\mathbb{R}^n)}dt\right)\\
&\leq & \frac{C}{\Gamma(k-s_1/2)}\left( T^{\frac{s-s_1}{2}} \mathcal{M} \left((-\Delta)^{\frac{s}{2}}(f)\right)(x)+ T^{\frac{-\beta-s_1}{2}}\|f\|_{\dot{B}^{-\beta, \infty}_{\infty}(\mathbb{R}^n)}\right).
\end{eqnarray*}
We fix now the parameter $T$ by the condition
$$T=\left(\frac{\|f\|_{\dot{B}^{-\beta, \infty}_{\infty}(\mathbb{R}^n)}}{\mathcal{M} \left((-\Delta)^{\frac{s}{2}} (f)\right)(x)}\right)^{ \frac{2}{\beta+s}},$$
and we obtain the following inequality
$$|(-\Delta)^{\frac{s_1}{2}}(f)(x)|\leq \frac{C}{\Gamma(k-s_1/2)}\mathcal{M} \left((-\Delta)^{\frac{s}{2}}(f)\right)^{1-\frac{s-s_1}{\beta+s}}(x)\|f\|^{\frac{s-s_1}{\beta+s}}_{\dot{B}^{-\beta, \infty}_{\infty}(\mathbb{R}^n)}.$$
Since $\frac{s-s_1}{\beta+s}=\theta$ we have
$$|(-\Delta)^{\frac{s_1}{2}}(f)(x)|\leq \frac{C}{\Gamma(k-s_1/2)}\mathcal{M} \left((-\Delta)^{\frac{s}{2}}(f)\right)^{1-\theta}(x)\|f\|^{\theta}_{\dot{B}^{-\beta, \infty}_{\infty}(\mathbb{R}^n)},$$
from which we deduce
$$|(-\Delta)^{\frac{s_1}{2}}(f)(x)|^{\frac{1}{1-\theta}}\leq C\mathcal{M} \left((-\Delta)^{\frac{s}{2}} (f)\right)(x)\|f\|^{\frac{\theta}{1-\theta}}_{\dot{B}^{-\beta, \infty}_{\infty}(\mathbb{R}^n)}.$$
\hfill$\blacksquare$
\section{Proof of Theorem \ref{Theo2}}\label{Secc_ProofTheo2}
The starting point of our proof is given by the modified Hedberg inequality (\ref{HedbergNuevo1}) given in Theorem \ref{Theo1}:
$$|(-\Delta)^{\frac{s_1}{2}}(f)(x)|^{\frac{1}{1-\theta}}\leq C\mathcal{M} \left((-\Delta)^{\frac{s}{2}}(f)\right)(x)\|f\|^{\frac{\theta}{1-\theta}}_{\dot{B}^{-\beta, \infty}_{\infty}(\mathbb{R}^n)}.$$
Now we apply the norm $\|\cdot\|_{L^{p(\cdot)}}$ on both sides of the previous inequality to obtain
$$\left\||(-\Delta)^{\frac{s_1}{2}}(f)(x)|^{\frac{1}{1-\theta}}\right\|_{L^{p(\cdot)}(\mathbb{R}^n)}\leq C\|f\|^{\frac{\theta}{1-\theta}}_{\dot{B}^{-\beta, \infty}_{\infty}(\mathbb{R}^n)} \|\mathcal{M} \left((-\Delta)^{\frac{s}{2}} (f)\right)\|_{L^{p(\cdot)}(\mathbb{R}^n)}.$$
Note that by the identity (\ref{PuissanceNorme}) we have $\||f|^\alpha\|_{L^{p(\cdot)}(\mathbb{R}^n)}=\|f\|_{L^{\alpha p(\cdot)}(\mathbb{R}^n)}^\alpha$ with $\frac{1}{ p^-}\leq \alpha$,  thus since $\frac{1}{p^-}\leq \frac{1}{1-\theta}$ (recall that $0<\theta<1$ and $1<p^-$), we can write 
$$\left\|(-\Delta)^{\frac{s_1}{2}}(f)\right\|_{L^{p(\cdot)/(1-\theta)}(\mathbb{R}^n)}^{\frac{1}{1-\theta}}\leq C\|f\|^{\frac{\theta}{1-\theta}}_{\dot{B}^{-\beta, \infty}_{\infty}(\mathbb{R}^n)} \|\mathcal{M} \left((-\Delta)^{\frac{s}{2}} (f)\right)\|_{L^{p(\cdot)}(\mathbb{R}^n)},$$
and since by hypothesis we have $q(\cdot)=p(\cdot)/(1-\theta)$, we obtain
$$\left\|(-\Delta)^{\frac{s_1}{2}}(f)\right\|_{L^{q(\cdot)}(\mathbb{R}^n)}^{\frac{1}{1-\theta}}\leq C\|f\|^{\frac{\theta}{1-\theta}}_{\dot{B}^{-\beta, \infty}_{\infty}(\mathbb{R}^n)} \|\mathcal{M}\left((-\Delta)^{\frac{s}{2}} (f)\right)\|_{L^{p(\cdot)}(\mathbb{R}^n)},$$
but since $p\in \mathcal{P}^{log}(\mathbb{R}^n)$, then the maximal function $\mathcal{M}$ is bounded in $L^{p(\cdot)}(\mathbb{R}^n)$, and taking the ($1-\theta$)-power of the previous inequality, we have
$$\left\|(-\Delta)^{\frac{s_1}{2}}(f)\right\|_{L^{q(\cdot)}(\mathbb{R}^n)}\leq C'\|f\|^{\theta}_{\dot{B}^{-\beta, \infty}_{\infty}(\mathbb{R}^n)} \|(-\Delta)^{\frac{s}{2}}(f)\|_{L^{p(\cdot)}(\mathbb{R}^n)}^{1-\theta},$$
which is the desired estimate. \hfill $\blacksquare$
\section{Proofs of the Theorems \ref{Theo3} and \ref{Theo4}}\label{Secc_ProofTheo3}
The proof of this theorem is based on the conclusion of Theorem \ref{Theo2} and on several embeddings between spaces of variable exponents. Indeed, since by hypothesis we have $f\in \dot{\mathcal{S}}^{s, p(\cdot)}_\mathfrak{p}(\mathbb{R}^n)$ (which is defined in (\ref{MixedSobolev})) then, in one hand, we have $f\in  \dot{W}^{s,\mathfrak{p}}(\mathbb{R}^n)$ with $0<s<n/\mathfrak{p}$, and due to the usual Sobolev embeddings we get
$$\dot{W}^{s,\mathfrak{p}}(\mathbb{R}^n)\subset L^r(\mathbb{R}^n)\qquad \mbox{with } \frac{n}{r}=\frac{n}{\mathfrak{p}}-s,$$
moreover, we also have the space inclusion 
$$L^r(\mathbb{R}^n)\subset \dot{B}^{-\beta,\infty}_{\infty}(\mathbb{R}^n)\qquad \mbox{with } \beta=\frac{n}{r},$$
(see \cite{Triebel} for a proof of this fact, see also \cite[Chapter 2]{Bahouri} for more details), from which we obviously get 
$$\dot{W}^{s,\mathfrak{p}}(\mathbb{R}^n)\subset \dot{B}^{-\beta,\infty}_{\infty}(\mathbb{R}^n),$$ 
where $\beta=\frac{n}{\mathfrak{p}}-s$ and thus we have $f\in \dot{B}^{-\beta,\infty}_{\infty}(\mathbb{R}^n)$.\\

Now, in the other hand, since we have $f\in \dot{W}^{s,p(\cdot)}(\mathbb{R}^n)$ where $0<s<n/p^+$ and where $p\in \mathcal{P}^{\log}(\mathbb{R}^n)$ with $1<p^-\leq p^+<+\infty$, we can apply Theorem \ref{Theo2} with $s_1=0$, $\beta=\frac{n}{\mathfrak{p}}-s$ and $\theta=\frac{s}{\beta+s}=\frac{s\mathfrak{p}}{n}$, to the function $f\in \dot{W}^{s,p(\cdot)}(\mathbb{R}^n)\cap \dot{B}^{-\beta,\infty}_{\infty}(\mathbb{R}^n)$ to obtain
\begin{eqnarray}
\|f\|_{L^{p(\cdot)/(1-\theta)}(\mathbb{R}^n)}&\leq &C \|f\|_{\dot{W}^{s,p(\cdot)}(\mathbb{R}^n)}^{1-\theta}\|f\|_{\dot{B}^{-\beta,\infty}_{\infty}(\mathbb{R}^n)}^\theta\notag \\
&\leq &C \|f\|_{\dot{W}^{s,p(\cdot)}(\mathbb{R}^n)}^{1-\theta}\|f\|_{\dot{W}^{s,\mathfrak{p}}(\mathbb{R}^n)}^{\theta},\label{InegaliteInterpol}
\end{eqnarray}
where in the last line we used the embedding $\dot{W}^{s,\mathfrak{p}}(\mathbb{R}^n)\subset \dot{B}^{-\beta,\infty}_{\infty}(\mathbb{R}^n)$ obtained previously. \\

We define now $\sigma(\cdot)=\frac{p(\cdot)}{1-\theta}$ but since $1-\theta=\frac{n-s\mathfrak{p}}{n}$, we have $\sigma(\cdot)=\frac{np(\cdot)}{n-s\mathfrak{p}}$ and recalling the definition of the quantity $\|\cdot\|_{\dot{\mathcal{S}}^{s, p(\cdot)}_\mathfrak{p}(\mathbb{R}^n)}$ given in (\ref{MixedSobolev}) and we easily obtain
$$\|f\|_{L^{\sigma(\cdot)}(\mathbb{R}^n)}\leq C \|f\|_{\dot{\mathcal{S}}^{s, p(\cdot)}_\mathfrak{p}(\mathbb{R}^n)},$$
which ends the proof of Theorem \ref{Theo3}. \hfill $\blacksquare$\\

\noindent{\bf Proof of Theorem \ref{Theo4}.} From the estimate (\ref{InegaliteInterpol}), the Theorem \ref{Theo4} follows easily: indeed, replacing formally $f$ by $I_s(f)$ we have
\begin{equation}\label{InegaliteInterpol1}
\|I_s(f)\|_{L^{\sigma(\cdot)}(\mathbb{R}^n)}\leq C \|f\|_{L^{p(\cdot)}(\mathbb{R}^n)}^{1-\theta}\|f\|_{L^{\mathfrak{p}}(\mathbb{R}^n)}^\theta,
\end{equation}
and by the definition of the quantity $\|\cdot\|_{\mathcal{L}^{p(\cdot)}_\mathfrak{p}(\mathbb{R}^n)}$ given in (\ref{MixedLebesgue}) we obtain the wished estimate (\ref{Inegalite3}) and Theorem \ref{Theo4} is proven. \hfill $\blacksquare$
\begin{Remarque}
Observe that the previous estimate (\ref{InegaliteInterpol1}) can also be obtained by interpolation. 
\end{Remarque}
Indeed, from the hypotheses of Theorem \ref{Theo4}, \emph{i.e.} $f \in L^{p(\cdot)}(\mathbb{R}^n)$ and $f\in L^{\mathfrak{p}}(\mathbb{R}^n)$ we can deduce from the usual Hardy-Littlewood inequalities (\ref{PotentialRieszVariable}) and (\ref{Sobolev3}) that 
$$\|I_s(f)\|_{L^{q(\cdot)}(\mathbb{R}^n)}\leq C \|f\|_{L^{p(\cdot)}(\mathbb{R}^n)}\qquad  \mbox{and} \qquad\|I_s(f)\|_{L^{\mathfrak{q}}(\mathbb{R}^n)}\leq C\|f\|_{L^{\mathfrak{p}}(\mathbb{R}^n)},$$
where $q(\cdot)=\frac{np(\cdot)}{n-sp(\cdot)}$ and $\mathfrak{q}=\frac{n \mathfrak{p}}{n-s\mathfrak{p}}$, thus, applying the complex interpolation theory (see Theorem 3.47 of \cite{Cruz1}) we can recover the space $L^{\sigma(\cdot)}(\mathbb{R}^n)$ and the estimate (\ref{InegaliteInterpol1}) above. 
\section{Orlicz spaces}\label{Secc_Orlicz}
We consider here another generalization of the usual Lebesgue spaces where Sobolev inequalities have been studied. Indeed, in the framework of \emph{Orlicz spaces}, some variants of the Sobolev inequalities (\ref{Sobolev1}) and the Hardy-Littlewood-Sobolev inequalities (\ref{Sobolev3}) are available (see \emph{e.g.} \cite{Cianchi}, \cite{Nakai} and \cite{Derigoz}) and their proofs also rely in suitable versions of the Hedberg inequality (\ref{Hedberg}).\\

 Let us recall that if $a:[0,+\infty[\longrightarrow[0,+\infty[$ is a left-continuous non decreasing function with $a(0)=0$, we can consider the corresponding \emph{Young function} $A(t)=\displaystyle{\int_{0}^ta(\sigma)d\sigma}$ and then the Orlicz space $L^A(\mathbb{R}^n)$ associated to the function $A$ is defined as the set of measurable functions $f:\mathbb{R}^n\longrightarrow \mathbb{R}$ such that the Luxemburg norm
\begin{equation*}
\|f\|_{L^A(\mathbb{R}^n)}=\inf\left\{\lambda > 0: \, \int_{\mathbb{R}^n}A(|f(x)|/\lambda)dx\leq1\right\},
\end{equation*}
is finite. The previous expression is of course very similar to (\ref{Def_LuxNormLebesgue}) and as we can easily see here that if $A(t)=t^p$ for $1\leq p<+\infty$, we recover the classical Lebesgue spaces. Since the quantity $\|\cdot\|_{L^A(\mathbb{R}^n)}$ is a norm, we can expect some usual properties: for exemple, just as in (\ref{ConstanteNormeLuxem1}), for a constant $\lambda\in \mathbb{R}$ we have 
$$\|\lambda f\|_{L^A(\mathbb{R}^n)}=|\lambda|\; \|f\|_{L^A(\mathbb{R}^n)},$$
and if $f,g$ are two measurable functions such that $|f|\leq |g|$ a.e., then we have the same order-reserving property given in (\ref{OrdreNormeLuxem1})
$$\|f\|_{L^A(\mathbb{R}^n)}\leq \|g\|_{L^A(\mathbb{R}^n)}.$$
However, the property (\ref{PuissanceNorme}) for the Lebesgue spaces of variable exponent should be handled more carefully and for this we will use the following rescaling property as defined in Section 3 of \cite{RaSam}: for any real $\sigma>0$, we define the space $L^A_\sigma(\mathbb{R}^n)$ by the condition
$$L^A_\sigma(\mathbb{R}^n)=\{f:\mathbb{R}^n\longrightarrow \mathbb{R}: \|f\|_{L^A_\sigma(\mathbb{R}^n)}<+\infty\},$$
where 
\begin{equation}\label{Rescal}
\|f\|_{L^A_\sigma(\mathbb{R}^n)}=\inf\left\{\lambda > 0: \, \int_{\mathbb{R}^n}A_\sigma(|f(x)|/\lambda)dx\leq1\right\},
\end{equation}
with $A_\sigma(t)=A(t^\sigma)$. With this definition of the functional $\|\cdot\|_{L^A_\sigma(\mathbb{R}^n)}$ we have the following property:
\begin{equation}\label{Rescal1}
\||f|^\sigma\|_{L^A(\mathbb{R}^n)}=\|f\|_{L^A_\sigma(\mathbb{R}^n)}^\sigma.
\end{equation}
See Lemma 3.2 of \cite{RaSam} for a proof of this fact.\\

Again, and just as for the Lebesgue spaces of variable exponent considered before, most of the usual tools and inequalities are harder to use in the setting of Orlicz spaces than in the classical framework as they strongly depend on the properties of the Young function $A$. \\

In order to study Sobolev inequalities, the approach displaying some adapted versions of the Hedberg inequality is commonly used and this requires two ingredients: the boundedness of the maximal functions and some point-wise estimate. For the first ingredient, it is classical to impose the following restrictions over the Young function $A$: a Young function $A$ is said to satisfy the $\nabla_2$-condition, denoted also by $A\in \nabla_2$, if
$$A(r)\leq\frac{1}{2C} A(Cr),\qquad r\geq 0,$$
for some $C > 1$, then if $A\in  \nabla_2$ we have
\begin{equation*}
\|\mathcal{M}(f)\|_{L^{A}(\mathbb{R}^n)}\leq C\|f\|_{L^{A}(\mathbb{R}^n)},
\end{equation*}
see \cite{Cianchi0} for a proof of this fact, see also \cite[Theorem 2]{Derigoz} and the reference there in for more details on the boundedness of the maximal functions in this setting.\\ 

Once we have at our disposal this boundedness property for the maximal functions, the Sobolev inequalities can be studied via suitable versions of the Hedberg inequality. Indeed, in \cite{Cianchi} the following pointwise Hedberg-type estimate is proven
\begin{equation}\label{Cianchi}
I_s(f)(x)\leq C\|f\|_{L^A(\mathbb{R}^n)}H_s\left(\frac{\mathcal{M}(f)(x)}{\|f\|_{L^A(\mathbb{R}^n)}}\right), \qquad \mbox{for } 0<s<n, 
\end{equation}
where the function $H_s$ is defined by the formula
$$H_s(\tau)=\left(\int_0^\tau\left(\frac{r}{A(r)}\right)dr\right)^{(n-s)/n},$$
thus, setting $A_s(\tau)=A(H_s(\tau)^{-1})$ for $\tau\geq 0$, the following version of the Hardy-Littlewood-Sobolev inequalites are obtained:
\begin{equation*}
\|I_s(f)\|_{L^{A_s}(\mathbb{R}^n)}\leq C\|f\|_{L^A(\mathbb{R}^n)}.
\end{equation*}
Another method to obtain Hedberg-type inequalities is displayed in Theorem 7.1 of \cite{Nakai}, where some inequalities involving fractional Riesz potentials are obtained. \\

Let us remark now that, contrary to the inequality (\ref{Cianchi}), the modified Hedberg inequality (\ref{HedbergNuevo1}) proposed in this article does not require the presence of any suitable Young function and this special structure allows us to obtain in a very straightforward manner the following inequality:
\begin{Theoreme}[Sobolev-like inequalities for Orlicz spaces]\label{Theo5}
Let $A$ be a Young function such that $A\in \nabla_2$. Consider $f:\mathbb{R}^n\longrightarrow \mathbb{R}$ a measurable function such that $(-\Delta)^{\frac s2}(f)\in L^{A}(\mathbb{R}^n)$ and $f\in \dot{B}^{-\beta,\infty}_{\infty}(\mathbb{R}^n)$ for some $\beta>0$. Define $\theta=\frac{s-s_1}{\beta+s}<1$ for $0\leq s_1<s$. Then we have $(-\Delta)^{\frac{s_1}{2}}(f)\in L^{A}_{(1-\theta)}(\mathbb{R}^n)$ where the space $L^{A}_{(1-\theta)}(\mathbb{R}^n)$ is defined as in (\ref{Rescal})-(\ref{Rescal1}). Moreover, the following inequality holds true
\begin{equation}\label{Inegalite1Orlicz}
\|(-\Delta)^{\frac{s_1}{2}}(f)\|_{L^{A}_{(1-\theta)}(\mathbb{R}^n)}\leq C \|(-\Delta)^{\frac{s}{2}}(f)\|_{L^{A}(\mathbb{R}^n)}^{1-\theta}\|f\|_{\dot{B}^{-\beta,\infty}_{\infty}(\mathbb{R}^n)}^\theta.
\end{equation}
\end{Theoreme}
{\bf Proof.} The arguments follow very closely those given in Section \ref{Secc_ProofTheo2} and we given them for the sake of completeness: from inequality (\ref{HedbergNuevo1}) we have
$$|(-\Delta)^{\frac{s_1}{2}}(f)(x)|^{\frac{1}{1-\theta}}\leq C\mathcal{M} \left((-\Delta)^{\frac{s}{2}}(f)\right)(x)\|f\|^{\frac{\theta}{1-\theta}}_{\dot{B}^{-\beta, \infty}_{\infty}(\mathbb{R}^n)},$$
and we apply the Luxemburg norm $\|\cdot\|_{L^{A}}$ associated to the Young function $A$ on both sides of the inequality to get
$$\left\||(-\Delta)^{\frac{s_1}{2}}(f)|^{\frac{1}{1-\theta}}\right\|_{L^{A}(\mathbb{R}^n)}\leq C\|f\|^{\frac{\theta}{1-\theta}}_{\dot{B}^{-\beta, \infty}_{\infty}(\mathbb{R}^n)} \|\mathcal{M}\left((-\Delta)^{\frac{s}{2}} (f)\right)\|_{L^{A}(\mathbb{R}^n)}.$$
We use now the rescaling property (\ref{Rescal1}) to obtain the identity
$$\||(-\Delta)^{\frac{s_1}{2}}(f)|^{\frac{1}{1-\theta}}\|_{L^A(\mathbb{R}^n)}=\|(-\Delta)^{\frac{s_1}{2}}(f)\|_{L^A_{(1-\theta)}(\mathbb{R}^n)}^{\frac{1}{1-\theta}},$$
and since $A\in \nabla_2$, the Hardy-Littlewood maximal operator is bounded in the space $L^A(\mathbb{R}^n)$, we can write 
$$\|(-\Delta)^{\frac{s_1}{2}}(f)\|_{L^A_{(1-\theta)}(\mathbb{R}^n)}^{\frac{1}{1-\theta}}\leq C\|f\|^{\frac{\theta}{1-\theta}}_{\dot{B}^{-\beta, \infty}_{\infty}(\mathbb{R}^n)} \|(-\Delta)^{\frac{s}{2}} (f)\|_{L^{A}(\mathbb{R}^n)},$$
which is the desired estimate (\ref{Inegalite1Orlicz}). \hfill $\blacksquare$
\begin{Remarque}
Of course, if $s_1=0$ we have the estimate
$$\|f\|_{L^A_{(1-\theta)}(\mathbb{R}^n)}\leq C\|(-\Delta)^{\frac{s}{2}} (f)\|_{L^{A}(\mathbb{R}^n)}^{1-\theta}\|f\|^{\theta}_{\dot{B}^{-\beta, \infty}_{\infty}(\mathbb{R}^n)},$$
which constitutes a new variant of the usual Sobolev inequalities in the framework of Orlicz spaces.
\end{Remarque}
To end this section, note that in Lemma 9 of \cite{Derigoz} a variant of the inequality (\ref{Cianchi}) for generalized fractional integral operators is obtained where the function $H_s$ is replaced by a suitable Orlicz-Morrey space of the third kind (see Section 5-6 of \cite{Derigoz}). This result is in some sense close to ours (the function $H_s$ is replaced by another functional space) but the technics displayed in the mentioned article are very different from our approach. 
\section{Some possible generalizations}\label{Secc_General}
The modified Hedberg inequality (\ref{HedbergNuevo1}) and the subsequent Sobolev-like estimates can be applied in several settings. In this article we studied Lebesgue spaces with variable exponents and Orlicz spaces over $\mathbb{R}^n$ but we also can consider the following frameworks:
\begin{itemize}
\item[$\bullet$] Instead of the usual base space $\mathbb{R}^n$, it is possible to consider stratified Lie groups $\mathbb{G}$ such as the Heisenberg group (see \cite{Stein2} for the details).
\item[$\bullet$] As long as the boundedness of the maximal function is preserved, several type of weights can be considered in the previous inequalities. See \cite[Chapter 5]{Diening} for a generalized Muckenhoupt condition in the setting of Lebesgue spaces of variable exponents. See \cite{Krbec} for a theory of weighted Orlicz spaces.
\end{itemize}

\appendix
\mysection{Appendix}\label{AppendixA}

An alternative proof of Theorem \ref{Theo2} relies in the characterization of the Lebesgue and Sobolev spaces of variable exponents using the Littlewood-Paley decomposition. In the classical setting, the representation of these spaces using dyadic blocs is well known (see the books \cite{Bahouri} and \cite{Grafakos}). For spaces of variable exponents, this theory is given in \cite{Almeida}, see also Chapter 12 of the book \cite{Diening}. Let us point out that the identification of the spaces given by a Littlewood-Paley decomposition with the spaces used here (in the non-homogeneous case) is done in the article \cite{Diening1}.\\

Let us briefly recall the Littlewood-Paley decomposition: consider $\varphi\in \mathcal{S}(\mathbb{R}^n, \mathbb{R})$ such that 
$\widehat{\varphi}(\xi)=1$ if $|\xi|\leq 1/2$ and $\widehat{\varphi}(\xi)=0$ if $|\xi|>1$ and for $j\in \mathbb{Z}$ define the function $\varphi_j$ by the expression $\varphi_j(x)=2^{-jn}\varphi(2^{-j}x)$. Consider the function $\psi$ which is given by the formula $\widehat{\psi}(\xi)=\widehat{\varphi}(\xi/2)-\widehat{\varphi}(\xi)$ and for all $j\in \mathbb{Z}$ we set $\widehat{\psi}_j(\xi)=\widehat{\psi}(2^{-j}\xi)=\widehat{\varphi}_{j+1}(\xi)-\widehat{\varphi}_{j}(\xi)$, we have then 
$$\sum_{j\in \mathbb{Z}}\widehat{\psi}_j(\xi)\equiv 1\quad \mbox{for all } \xi\neq 0.$$
Now, for all $j\in \mathbb{Z}$, the dyadic-bloc operator $\Delta_j$ is defined by the formula $\Delta_j(f)=f\ast \psi_j$ and we have the formula
$$f=\sum_{j\in \mathbb{Z}}\Delta_j(f),$$
where the convergence of the sum must be considered in $\mathcal{S}'(\mathbb{R}^n)$ modulo the polynomials $\mathbb{C}[X]$.\\

Now for $p \in \mathcal{P}(\mathbb{R}^n)$ such that $1<p^-\leq p^+<+\infty$ and for $0<s<+\infty$, we have the following characterizations of variable exponents spaces: for all $f\in \mathcal{S}'/\mathbb{C}[X]$ we have
$$\|f\|_{L^{p(\cdot)}(\mathbb{R}^n)}\simeq \mbox{\footnotesize$\left\|\left(\displaystyle{\sum_{j\in \mathbb{Z}}}|\Delta_{j}(f)(x)|^{2}\right)^{\frac{1}{2}}\right\|_{L^{p(\cdot)}(\mathbb{R}^n)}$}\mbox{and} \quad \|f\|_{\dot{\mathcal{W}}^{s,p(\cdot)}(\mathbb{R}^n)}\simeq\mbox{\footnotesize$\left\|\left(\displaystyle{\sum_{j\in \mathbb{Z}}}2^{2sj}|\Delta_{j}(f)(x)|^{2}\right)^{\frac{1}{2}}\right\|_{L^{p(\cdot)}(\mathbb{R}^n)}$}.$$
Note that the spaces $\dot{\mathcal{W}}^{s,p(\cdot)}(\mathbb{R}^n)$ given above are defined modulo the polynomials and thus they are \emph{not} equivalent to the spaces $\dot{W}^{s,p(\cdot)}(\mathbb{R}^n)$ given in (\ref{Def_Sobolev}). However, in the framework of the Theorem \ref{Theo2}, we are considering functions that also belong to the Besov space 
$\dot{B}^{-\beta,\infty}_{\infty}(\mathbb{R}^n)$, which can be characterized by the equivalent quantity 
\begin{equation}\label{Def_besov_LP}
\|f\|_{\dot{B}^{-\beta,\infty}_{\infty}(\mathbb{R}^n)}\simeq\underset{j\in\mathbb{Z}}{\sup\;} 2^{-\beta j}\|\Delta_{j}(f)\|_{L^\infty(\mathbb{R}^n)},
\end{equation}
and this fact allows us to avoid this unpleasant issue related to the different definitions of homogeneous spaces as we have the equivalence of spaces\footnote{Polynomials are excluded here and this fact can be easily seen from the characterization (\ref{Def_BesovThermic}) of Besov spaces.}
$\dot{W}^{s,p(\cdot)}\cap \dot{B}^{-\beta,\infty}_{\infty}\simeq \dot{\mathcal{W}}^{s,p(\cdot)}\cap \dot{B}^{-\beta,\infty}_{\infty}$.
\quad\\

With this short introduction, we can proof Theorem \ref{Theo2} using the tools related to the Littlewood-Paley decomposition. We start with the following interpolation result 
\begin{Lemme} Let $(a_{j})_{j\in \mathbb{Z}}$ be a sequence and set $s=(1-\theta) s_0+ \theta s_{1}$ with $0<\theta< 1$ and $s_0\neq s_{1}$.  Then for all $r,r_{1},r_{2}\in [1,+\infty]$ we have the interpolation estimate:
\begin{equation*}
\|2^{js}a_{j}\|_{\ell^{r}}\leq C\|2^{js_{0}}a_{j}\|_{\ell^{r_{1}}}^{1-\theta}\|2^{js_{1}}a_{j}\|_{\ell^{r_{2}}}^{\theta}.
\end{equation*}
\end{Lemme}
See \cite{Bergh} for a proof of this interpolation inequality.\\

If we apply this lemma to the dyadic blocs $\Delta_{j}(f)$ with $s_1=(1-\theta)s+\theta (-\beta)$ and $r=r_{1}=2$ and $r_{2}=+\infty$, we obtain
\begin{eqnarray*}
\left(\sum_{j\in \mathbb{Z}}2^{2s_1j}|\Delta_{j}(f)(x)|^{2}\right)^{\frac{1}{2(1-\theta)}}&\leq &C \left(\sum_{j\in \mathbb{Z}}2^{2sj}|\Delta_{j}(f)(x)|^{2}\right)^{\frac{1}{2}}
\left(\underset{j\in\mathbb{Z}}{\sup\;} 2^{-\beta j}|\Delta_{j}(f)(x)|\right)^{\frac{\theta}{1-\theta}}\\
&\leq & C \left(\sum_{j\in \mathbb{Z}}2^{2sj}|\Delta_{j}(f)(x)|^{2}\right)^{\frac{1}{2}}\|f\|_{\dot{B}^{-\beta,\infty}_\infty(\mathbb{R}^n)}^{\frac{\theta}{1-\theta}},
\end{eqnarray*}
where in the last line we used the characterization of the Besov space $\dot{B}^{-\beta,\infty}_\infty(\mathbb{R}^n)$ in terms of the dyadic blocs given in (\ref{Def_besov_LP}). Now, we take the Luxemburg $L^{p(\cdot)}$-norm to get
$$\left\| \left(\sum_{j\in \mathbb{Z}}2^{2s_1j}|\Delta_{j}(f)|^{2}\right)^{\frac12 \frac{1}{1-\theta}}\right\|_{L^{p(\cdot)}(\mathbb{R}^n)}\leq C \left\| \left(\sum_{j\in \mathbb{Z}}2^{2sj}|\Delta_{j}(f)|^{2}\right)^{\frac12}
\right\|_{L^{p(\cdot)}(\mathbb{R}^n)}\|f\|_{\dot{B}^{-\beta,\infty}_\infty(\mathbb{R}^n)}^{\frac{\theta}{1-\theta}}.$$
Then, using the property (\ref{PuissanceNorme}) we obtain 
$$\left\| \left(\sum_{j\in \mathbb{Z}}2^{2s_1j}|\Delta_{j}(f)|^{2}\right)^{\frac12 }\right\|_{L^{\frac{p(\cdot)}{1-\theta}}(\mathbb{R}^n)}^{\frac{1}{1-\theta}}\leq C \left\| \left(\sum_{j\in \mathbb{Z}}2^{2sj}|\Delta_{j}(f)|^{2}\right)^{\frac12}
\right\|_{L^{p(\cdot)}(\mathbb{R}^n)}\|f\|_{\dot{B}^{-\beta,\infty}_\infty(\mathbb{R}^n)}^{\frac{\theta}{1-\theta}}.$$
To finish, we recall that $q(\cdot)=\frac{p(\cdot)}{1-\theta}$ and using the characterization of Sobolev spaces via the Littlewood-Paley theory we can write
$$\|f\|_{\dot{\mathcal{W}}^{s_1,q(\cdot)}(\mathbb{R}^n)}\leq C \|f\|_{\dot{\mathcal{W}}^{s,p(\cdot)}(\mathbb{R}^n)}^{1-\theta}\|f\|_{\dot{B}^{-\beta,\infty}_{\infty}(\mathbb{R}^n)}^\theta.$$
\hfill $\blacksquare$
\begin{Remarque}
Let us note here that this second proof of the Sobolev-like estimates relies in the Littlewood-Paley theory which is available in the setting of Lebesgue spaces of variable exponent over the euclidean space $\mathbb{R}^n$. But this is not always the case if we consider general spaces over other spaces than $\mathbb{R}^n$. In this sense the first proof based in the modified Hedberg inequality (\ref{HedbergNuevo1}) seems more robust and simple to display.
\end{Remarque}


\begin{thebibliography}{2}
\bibitem{Almeida}
A. \textsc{Almeida}, L. \textsc{Diening}, P. \textsc{H\"ast\"o}. \emph{Homogeneous variable exponent Besov and Triebel–Lizorkin spaces}. Mathematische Nachrichten. 291:1177–1190 (2018).
\bibitem{Bahouri}
H. \textsc{Bahouri}, J.-Y. \textsc{Chemin} and R. \textsc{Danchin}. \emph{Fourier Analysis and Nonlinear Partial Differential Equations}, Grundlehren der mathematischen Wissenschaften 343, Springer (2011).
\bibitem{Bergh}
J. \textsc{Bergh} and J. \textsc{Lofstrom}. \emph{Interpolation Spaces}, Grundlehren der mathematischen Wissenschaften 223, Springer (1976).
\bibitem{Brezis}
H. \textsc{Brezis} and P. \textsc{Mironescu}. \emph{Gagliardo-Nirenberg inequalities and non-inequalities: The full story}.  Annales de l'Institut Henri Poincaré - Non Linear Analysis 35, 1355-1376.  (2018).
\bibitem{Capone}
C. \textsc{Capone}, D. \textsc{Cruz-Uribe}, A. \textsc{Fiorenza}. \emph{The fractional maximal operator and fractional integrals on variable Lp spaces}. Rev. Mat. Iberoamericana 23, no. 3, 743–770 (2007).
\bibitem{Cianchi0}
A. \textsc{Cianchi}. \emph{Strong and weakly inequalities for some classical operators in Orlicz spaces}. J. Lond. Math. Soc. 60(1), 187–202 (1999).
\bibitem{Cianchi}
A. \textsc{Cianchi}. \emph{Optimal Orlicz-Sobolev embeddings}. Rev. Mat. Iberoamericana 20, no. 2, 427--474 (2004).
\bibitem{Cruz1}
D. \textsc{Cruz-Uribe}, A. \textsc{Fiorenza}. \emph{Variable Lebesgue Spaces}. Birkh\"auser (2013).
\bibitem{Derigoz}
F. \textsc{Deringoz} \emph{et al.}. \emph{Generalized fractional maximal and integral operators on Orlicz and generalized Orlicz–Morrey spaces of the third kind}. Positivity 23:727–757 (2019).
\bibitem{Diening}
L. \textsc{Diening}, P. \textsc{Harjulehto},  P. \textsc{H\"ast\"o},  M. \textsc{Ruzicka.} \emph{Lebesgue and Sobolev spaces with variable exponents}. Lecture Notes in Mathematics, 2017, Springer (2011).
\bibitem{Diening1}
L. \textsc{Diening}, P. \textsc{H\"ast\"o}, and S. \textsc{Roudenko}. \emph{Function spaces of variable smoothness and integrability}. J. Funct. Anal., 256(6):1731–1768, (2009).
\bibitem{GMO}
P. \textsc{G\'erard}, Y. \textsc{Meyer} \& F. \textsc{Oru}. \emph{In\'egalit\'es de Sobolev Pr\'ecis\'ees}. S\' eminaire sur les Equations aux
D\'eriv\'ees Partielles, 1996-1997, Exp. No. IV, \' Ecole Polytech., Palaiseau.
\bibitem{Grafakos}
L. \textsc{Grafakos}. \emph{Classical and Modern Fourier Analysis}. Prentice Hall (2004).
\bibitem{Hedberg}
L. \textsc{Hedberg}. \emph{On certain convolution inequalities}. Proc. Amer. Math. Soc. 36, 505–510 (1972).
\bibitem{Krbec}
V. \textsc{Kokilashvili} and M. Krbec. \emph{Weighted Inequalities in Lorentz and Orlicz Spaces}. World Scientific (1991).
\bibitem{Kolyada}
V. I. \textsc{Kolyada} and F. J. \textsc{P\'erez L\'azaro}.  \emph{On Gagliardo–Nirenberg Type Inequalities}. Journal of Fourier Analysis and Applications. volume 20, 577–607 (2014).
\bibitem{Nakai}
E. \textsc{Nakai}, \emph{Hardy-Littlewood maximal operator, singular integral operators and Riesz potentials on generalized Morrey spaces}. Math. Nachr. 166, 95--103 (1994).
\bibitem{RaSam}
H. \textsc{Rafeiro} and S. \textsc{Samko}. \emph{Maximal Operator with Rough Kernel in Variable Musielak–Morrey–Orlicz type Spaces, Variable Herz Spaces and Grand Variable Lebesgue Spaces}. Integr. Equ. Oper. Theory 89, 111–124 (2017).
\bibitem{Stein1}
E. M. \textsc{Stein}. \emph{Singular Integrals and Differentiability Properties of Functions}. Princeton Mathematical Series, 30. Princeton University Press (1970).
\bibitem{Stein2}
E. M. \textsc{Stein}. \emph{Harmonic Analysis}. Princeton University Press (1993).
\bibitem{Triebel}
H. \textsc{Triebel}. \emph{Theory of Function Spaces}. Birkh\"auser (1983).
\end{thebibliography}
\end{document}